\documentclass[12pt]{article}%
\usepackage{amsfonts}
\usepackage{sw20bams}
\usepackage{amsmath}
\usepackage{amssymb}
\usepackage{graphicx}%
\providecommand{\U}[1]{\protect\rule{.1in}{.1in}}
\begin{document}

\title{Second Best, Third Worst, Fourth in Line}
\author{Steven Finch}
\date{May 1, 2022}
\maketitle

\begin{abstract}
We investigate decomposable combinatorial labeled structures more fully,
focusing on the exp-log class of type $a=1$ or $1/2$. \ For instance, the
modal length of the second longest cycle in a random $n$-permutation is
$(0.2350...)n$, whereas the modal length of the second smallest component in a
random $n$-mapping is $2$ (conjecturally, given $n\geq434$). \ As in earlier
work, our approach is to establish how well existing theory matches
experimental data and to raise open questions.

\end{abstract}

\footnotetext{Copyright \copyright \ 2022 by Steven R. Finch. All rights
reserved.}Given a combinatorial object with $n$~nodes, our interest is in

\begin{itemize}
\item the size of its $r^{\text{th}}$ longest cycle or largest component,

\item the size of its $r^{\text{th}}$ shortest cycle or smallest component
\end{itemize}

\noindent where $r\geq2$. \ If the object has no $r^{\text{th}}$ component,
then its $r^{\text{th}}$ largest/smallest components are defined to have
length $0$. \ The case $r=1$ has attracted widespread attention
\cite{Fi1-tcs8, Fi2-tcs8}. \ Key to our prior study were recursive formulas
\cite{GG-tcs8, PR1-tcs8} for $L_{k,n}$ and $S_{k,n}$, the number of
$n$-objects whose largest and smallest components, respectively, have exactly
$k$ nodes, $1\leq k\leq n$. \ Different algorithms shall be used here. \ As
before, an $n$-object is chosen uniformly at random. \ For simplicity, we
discuss here only $n$-permutations and $n$-mappings (from $\{1,2,\ldots,n\}$
to $\{1,2,\ldots,n\}$). \ Let $c_{n}$ be the number of $n$-objects that are
connected, i.e., who possess exactly one component:%
\[
c_{n}=\left\{
\begin{array}
[c]{lll}%
(n-1)! &  & \text{for permutations,}\\
n!%
{\displaystyle\sum\limits_{j=1}^{n}}
\dfrac{n^{n-j-1}}{(n-j)!} &  & \text{for mappings.}%
\end{array}
\right.
\]
The total number of $n$-permutations and $n$-mappings is $n!$ and $n^{n}$,
respectively. \ For fixed $n$, the sequences $\{L_{k,n}:1\leq k\leq n\}$ and
$\{S_{k,n}:1\leq k\leq n\}$ constitute probability mass functions (upon
normalization) for $r=1$. \ Until recently, calculating analogous sequences
for $r\geq2$ seemed inaccessibly difficult.

The new algorithms, due to Heinz \cite{Hnz-tcs8}, accept as input the integer
$n$~and an ordered $r$-tuple $\ell$ of nonnegative integers, which may include
infinity. \ We write $\ell$ as a list $\{i_{1},i_{2},\ldots,i_{r}\}$. \ Given
a positive integer $j$, define$\ \ell^{j}$ to be the list obtained by

\begin{enumerate}
\item[(i)] appending $\ell$ with $j$,

\item[(ii)] sorting the $(r+1)$-tuple in ascendent order, and

\item[(iii)] removing its first element.
\end{enumerate}

\noindent Define $\ell_{j}$ in the same way as $\ell^{j}$ except for a revised
final step:

\begin{enumerate}
\item[(iii')] removing its last element.
\end{enumerate}

\noindent Note that the lengths of $\ell^{j}$ and $\ell_{j}$ are always equal
to the length of $\ell$. \ Let $p[n,\ell]$ and $q[n,\ell]$ denote row
polynomials in $x$ and $y$ associated with large and small components. The
algorithms are based on recursions%
\[
p[n,\ell]=\left\{
\begin{array}
[c]{lll}%
{\displaystyle\sum\limits_{j=1}^{n}}
c_{j}\,p[n-j,\ell^{j}]\dbinom{n-1}{j-1} & \bigskip & \text{if }n>0,\\
x^{i_{1}} &  & \text{if }n=0;
\end{array}
\right.
\]%
\[
q[n,\ell]=\left\{
\begin{array}
[c]{lll}%
{\displaystyle\sum\limits_{j=1}^{n}}
c_{j}\,q[n-j,\ell_{j}]\dbinom{n-1}{j-1} & \bigskip & \text{if }n>0,\\
y^{i_{r}} & \bigskip & \text{if }n=0\text{ and }i_{r}<\infty,\\
y^{0} &  & \text{if }n=0\text{ and }i_{r}=\infty.
\end{array}
\right.
\]
A computer algebra software package (e.g., Mathematica or Maple) makes exact
integer calculations for ample $n$ of $p[n,\ell]$ and $q[n,\ell]$ feasible.
\ These are demonstrated for $n=4$ in the next section, for the sake of concreteness.

Permutations belong to the exp-log class of type $a=1$, whereas mappings
belong to the exp-log class of type $a=1/2$. \ Explaining the significance of
the parameter $a>0$ would take us too far afield \cite{PR2-tcs8}. \ Let%
\[%
\begin{array}
[c]{ccc}%
E(x)=%
{\displaystyle\int\limits_{x}^{\infty}}
\dfrac{e^{-t}}{t}dt=-\operatorname{Ei}(-x), &  & x>0
\end{array}
\]
be the exponential integral. \ Define \cite{SL-tcs8, Shi-tcs8, FO-tcs8,
Gou-tcs8, ABT-tcs8, Pin-tcs8}%
\[
_{L}G_{a}(r,h)=\frac{\Gamma(a+1)a^{r-1}}{\Gamma(a+h)(r-1)!}%
{\displaystyle\int\limits_{0}^{\infty}}
x^{h-1}E(x)^{r-1}\exp\left[  -a\,E(x)-x\right]  dx,
\]%
\[
_{S}G_{a}(r,h)=\left\{
\begin{array}
[c]{lll}%
e^{-h\,\gamma}a^{r-1}/r! &  & \text{if }h=a,\\
\dfrac{\Gamma(a+1)}{(h-1)!(r-1)!}%
{\displaystyle\int\limits_{0}^{\infty}}
x^{h-1}\exp\left[  a\,E(x)-x\right]  dx &  & \text{if }h>a
\end{array}
\right.
\]
which are related to the $h^{\text{th}}$ moment of the $r^{\text{th}}$
largest/smallest component size (in this paper, rank $r=2,3$ or $4$; height
$h=1$ or $2$). \ While moment formulas are unerring for $L$, they are
\textbf{not} so for $S$. \ While $_{S}G_{a}$ is flawless for permutations (and
for what are called \textit{cyclations} \cite{Pip-tcs8}), a correction factor
$\sqrt{2}$ is needed for mappings.

For fixed $n$ and $r$, the coefficient sequences associated with polynomials%
\[%
\begin{array}
[c]{ccc}%
p[n,\{\underset{r}{\underbrace{0,0,\cdots,0}}\}], &  & 0\leq k\leq\left\lfloor
n/r\right\rfloor ;
\end{array}
\]%
\[%
\begin{array}
[c]{ccc}%
q[n,\{\underset{r}{\underbrace{\infty,\infty,\cdots,\infty}}\}], &  & 0\leq
k\leq n-r+1
\end{array}
\]
constitute probability mass functions (upon normalization). \ These have
corresponding means $_{L}\mu_{n,r}$, $_{S}\mu_{n,r}$ and variances $_{L}%
\sigma_{n,r}^{2}$, $_{S}\sigma_{n,r}^{2}$ given in the tables. \ We also
provide the median $_{L}\nu_{n,r}$ and mode $_{L}\vartheta_{n,r}$; evidently
$_{S}\nu_{n,r}$ and $_{S}\vartheta_{n,r}$ are bounded for permutations as
$n\rightarrow\infty$ (the trend of $_{S}\nu_{n,r}$ is less clear for
mappings). \ In table headings only, the following notation is used:
\[%
\begin{array}
[c]{ccccccccc}%
_{L}\widetilde{\mu}_{n,r}=\dfrac{_{L}\mu_{n,r}}{n}, &  & _{L}\widetilde
{\sigma}_{n,r}^{2}=\dfrac{_{L}\sigma_{n,r}^{2}}{n^{2}}, &  & _{L}%
\widetilde{\nu}_{n,r}=\dfrac{_{L}\nu_{n,r}}{n}, &  & _{L}\widetilde{\vartheta
}_{n,r}=\dfrac{_{L}\vartheta_{n,r}}{n}, &  & _{S}\widetilde{\nu}_{n,r}%
=\dfrac{_{S}\nu_{n,r}}{n},
\end{array}
\]%
\[%
\begin{array}
[c]{ccc}%
_{S}\widetilde{\mu}_{n,r}=\left\{
\begin{array}
[c]{lll}%
\dfrac{_{S}\mu_{n,r}}{\ln(n)^{r}} & \bigskip & \text{if }a=1,\\
\dfrac{_{S}\mu_{n,r}}{n^{1/2}\ln(n)^{r-1}} &  & \text{if }a=1/2,
\end{array}
\right.  &  & _{S}\widetilde{\sigma}_{n,r}^{2}=\left\{
\begin{array}
[c]{lll}%
\dfrac{_{S}\sigma_{n,r}^{2}}{n\ln(n)^{r-1}} & \bigskip & \text{if }a=1,\\
\dfrac{_{S}\sigma_{n,r}^{2}}{n^{3/2}\ln(n)^{r-1}} &  & \text{if }a=1/2.
\end{array}
\right.
\end{array}
\]
When $r=1$, the mode $_{L}\widetilde{\vartheta}_{n,1}$ is provably $1/2$ in
the limit as $n\rightarrow\infty$ for permutations (it is $1$ for mappings).
\ This limit is more interesting when $r=2$, as will soon be seen.

\section{Calculs \`{a} la Heinz}

As promised, we exhibit some hand calculations. \ It is easy to show directly
that $p[3,\{0,0\}]=2+4x$ for permutations and $17+10x$ for mappings (see
Section 3 of \cite{Fi1-tcs8}). \ More generally, $p[3,\{0,0\}]=c_{3}%
+c_{1}(c_{1}^{2}+3c_{2})x$. \ Let us compute $p[4,\{0,0\}]$ using Heinz's
algorithm. \ From%
\begin{align*}
p[2,\{1,1\}]  &  =c_{1}p[1,\{1,1\}]\tbinom{1}{0}+c_{2}p[0,\{1,2\}]\tbinom
{1}{1}\\
&  =c_{1}^{2}p[0,\{1,1\}]\tbinom{0}{0}+c_{2}x^{1}=\left(  c_{1}^{2}%
+c_{2}\right)  x,
\end{align*}%
\[
p[1,\{1,2\}]=c_{1}p[0,\{1,2\}]\tbinom{0}{0}=c_{1}x,
\]%
\[
p[0,\{1,3\}]=x
\]
we have%
\begin{align*}
p[3,\{0,1\}]  &  =c_{1}p[2,\{1,1\}]\tbinom{2}{0}+c_{2}p[1,\{1,2\}]\tbinom
{2}{1}+c_{3}p[0,\{1,3\}]\tbinom{2}{2}\\
&  =c_{1}\left(  c_{1}^{2}+c_{2}\right)  x+2c_{2}(c_{1}x)+c_{3}x=\left(
c_{1}^{3}+3c_{1}c_{2}+c_{3}\right)  x.
\end{align*}
Also, from%
\begin{align*}
p[2,\{0,2\}]  &  =c_{1}p[1,\{1,2\}]\tbinom{1}{0}+c_{2}p[0,\{2,2\}]\tbinom
{1}{1}\\
&  =c_{1}^{2}p[0,\{1,2\}]\tbinom{0}{0}+c_{2}x^{2}=c_{1}^{2}x+c_{2}x^{2},
\end{align*}%
\[
p[1,\{0,3\}]=c_{1}p[0,\{1,3\}]\tbinom{0}{0}=c_{1}x,
\]%
\[
p[0,\{0,4\}]=x^{0}=1
\]
we deduce%
\begin{align*}
p[4,\{0,0\}]  &  =c_{1}p[3,\{0,1\}]\tbinom{3}{0}+c_{2}p[2,\{0,2\}]\tbinom
{3}{1}+c_{3}p[1,\{0,3\}]\tbinom{3}{2}+c_{4}p[0,\{0,4\}]\tbinom{3}{3}\\
&  =c_{1}\left(  c_{1}^{3}+3c_{1}c_{2}+c_{3}\right)  x+3c_{2}(c_{1}^{2}%
x+c_{2}x^{2})+3c_{3}(c_{1}x)+c_{4}\\
&  =c_{4}+c_{1}\left(  c_{1}^{3}+6c_{1}c_{2}+4c_{3}\right)  x+3c_{2}^{2}%
x^{2}\\
&  =\left\{
\begin{array}
[c]{lll}%
6+15x+3x^{2} &  & \text{for permutations,}\\
142+87x+27x^{2} &  & \text{for mappings}%
\end{array}
\right.
\end{align*}
completing the argument.

It is likewise easy to show that $q[3,\{0,0\}]=2+y+3y^{2}$ for permutations
and $17+y+9y^{2}$ for mappings. \ More generally, $q[3,\{0,0\}]=c_{3}%
+c_{1}^{3}y+3c_{1}c_{2}y^{2}$. \ Let us compute $q[4,\{0,0\}]$ using Heinz's
algorithm. \ From%
\begin{align*}
q[2,\{1,1\}]  &  =c_{1}q[1,\{1,1\}]\tbinom{1}{0}+c_{2}q[0,\{1,1\}]\tbinom
{1}{1}\\
&  =c_{1}^{2}q[0,\{1,1\}]\tbinom{0}{0}+c_{2}y^{1}=\left(  c_{1}^{2}%
+c_{2}\right)  y,
\end{align*}%
\[
q[1,\{1,2\}]=c_{1}q[0,\{1,1\}]\tbinom{0}{0}=c_{1}y,
\]%
\[
q[0,\{1,3\}]=y^{3}%
\]
we have%
\begin{align*}
q[3,\{1,\infty\}]  &  =c_{1}q[2,\{1,1\}]\tbinom{2}{0}+c_{2}q[1,\{1,2\}]\tbinom
{2}{1}+c_{3}q[0,\{1,3\}]\tbinom{2}{2}\\
&  =c_{1}\left(  c_{1}^{2}+c_{2}\right)  y+2c_{2}(c_{1}y)+c_{3}y^{3}=\left(
c_{1}^{3}+3c_{1}c_{2}\right)  y+c_{3}y^{3}.
\end{align*}
Also, from%
\begin{align*}
q[2,\{2,\infty\}]  &  =c_{1}q[1,\{1,2\}]\tbinom{1}{0}+c_{2}q[0,\{2,2\}]\tbinom
{1}{1}\\
&  =c_{1}^{2}q[0,\{1,1\}]\tbinom{0}{0}+c_{2}y^{2}=c_{1}^{2}y+c_{2}y^{2},
\end{align*}%
\[
q[1,\{3,\infty\}]=c_{1}q[0,\{1,3\}]\tbinom{0}{0}=c_{1}y^{3},
\]%
\[
q[0,\{4,\infty\}]=y^{0}=1
\]
we deduce
\begin{align*}
q[4,\{\infty,\infty\}]  &  =c_{1}q[3,\{1,\infty\}]\tbinom{3}{0}+c_{2}%
q[2,\{2,\infty\}]\tbinom{3}{1}+c_{3}q[1,\{3,\infty\}]\tbinom{3}{2}%
+c_{4}q[0,\{4,\infty\}]\tbinom{3}{3}\\
&  =c_{1}\left(  \left(  c_{1}^{3}+3c_{1}c_{2}\right)  y+c_{3}y^{3}\right)
+3c_{2}(c_{1}^{2}y+c_{2}y^{2})+3c_{3}(c_{1}y^{3})+c_{4}\\
&  =c_{4}+c_{1}^{2}\left(  c_{1}^{2}+6c_{2}\right)  y+3c_{2}^{2}y^{2}%
+4c_{1}c_{3}y^{3}\\
&  =\left\{
\begin{array}
[c]{lll}%
6+7y+3y^{2}+8y^{3} &  & \text{for permutations,}\\
142+19y+27y^{2}+68y^{3} &  & \text{for mappings}%
\end{array}
\right.
\end{align*}
completing the argument.

\section{Modes \& Medians}

The mode of a continuous distribution is the location of its highest peak; the
median is its 50$^{\text{th}}$ percentile. \ The length $\Lambda_{r}$ of the
$r^{\text{th}}$ longest cycle in a random $n$-permutation has cumulative
probability%
\[
\lim_{n\rightarrow\infty}\mathbb{P}\left\{  \Lambda_{r}<x\cdot n\right\}
=\rho_{r}\left(  \frac{1}{x}\right)
\]
where $\rho_{r}(x)$ is the $r^{\text{th}}$ order Dickman function
\cite{KTP-tcs8}:%
\[%
\begin{array}
[c]{ccc}%
x\rho_{1}^{\prime}(x)+\rho_{1}(x-1)=0\text{ for }x>1, &  & \rho_{1}(x)=1\text{
for }0\leq x\leq1;
\end{array}
\]%
\[%
\begin{array}
[c]{ccc}%
x\rho_{r}^{\prime}(x)+\rho_{r}(x-1)=\rho_{r-1}(x-1)\text{ for }x>1, &  &
\rho_{r}(x)=1\text{ for }0\leq x\leq1
\end{array}
\]
and $r=2,3,4,\ldots$. \ For notational simplicity, let us write $\varphi
=\rho_{1}$ and $\psi=\rho_{2}$. \ Observe that $\rho_{r}$ should not be
confused with a different generalization $\rho_{a}$ discussed in
\cite{Fi1-tcs8, OPRW-tcs8}.

From%
\[%
\begin{array}
[c]{ccc}%
\varphi^{\prime}(x)=-\dfrac{\varphi(x-1)}{x}, &  & x>1
\end{array}
\]
we have%
\[%
\begin{array}
[c]{ccc}%
\varphi^{\prime}\left(  \dfrac{1}{x}\right)  =-\dfrac{\varphi\left(  \dfrac
{1}{x}-1\right)  }{\dfrac{1}{x}}, &  & 0<x<1
\end{array}
\]
hence the density $f(x)$ is%
\[
\dfrac{d}{dx}\varphi\left(  \dfrac{1}{x}\right)  =-x\,\varphi\left(  \dfrac
{1}{x}-1\right)  \left(  -\dfrac{1}{x^{2}}\right)  =\left\{
\begin{array}
[c]{lll}%
\dfrac{\varphi\left(  \dfrac{1}{x}-1\right)  }{x} & \bigskip & \text{if
}0<x\leq1/2,\\
\dfrac{1}{x} &  & \text{if }1/2<x<1.
\end{array}
\right.
\]
Also, from
\[%
\begin{array}
[c]{ccc}%
\varphi^{\prime\prime}(x)=\dfrac{\varphi(x-1)}{x^{2}}-\dfrac{\varphi^{\prime
}(x-1)}{x}=\dfrac{\varphi(x-1)}{x^{2}}+\dfrac{\varphi(x-2)}{x(x-1)}, &  & x>1
\end{array}
\]
we have%
\[%
\begin{array}
[c]{ccc}%
\varphi^{\prime\prime}\left(  \dfrac{1}{x}\right)  =\dfrac{\varphi\left(
\dfrac{1}{x}-1\right)  }{\dfrac{1}{x^{2}}}+\dfrac{\varphi\left(  \dfrac{1}%
{x}-2\right)  }{\dfrac{1}{x}\left(  \dfrac{1}{x}-1\right)  }, &  & 0<x<1
\end{array}
\]
hence (by the chain rule for second derivatives)%
\begin{align*}
\dfrac{d^{2}}{dx^{2}}\varphi\left(  \dfrac{1}{x}\right)   &  =\varphi^{\prime
}\left(  \dfrac{1}{x}\right)  \dfrac{2}{x^{3}}+\dfrac{1}{x^{4}}\varphi
^{\prime\prime}\left(  \dfrac{1}{x}\right) \\
&  =\dfrac{-2\varphi\left(  \dfrac{1}{x}-1\right)  }{x^{2}}+\frac
{\varphi\left(  \dfrac{1}{x}-1\right)  }{x^{2}}+\frac{\varphi\left(  \dfrac
{1}{x}-2\right)  }{x^{2}(1-x)}\\
&  =\left\{
\begin{array}
[c]{lll}%
\dfrac{1}{x^{2}(1-x)}-\dfrac{\varphi\left(  \dfrac{1}{x}-1\right)  }{x^{2}%
}>\dfrac{1}{x(1-x)}>0 & \bigskip & \text{if }1/3<x\leq1/2,\\
-\dfrac{1}{x^{2}}<0 &  & \text{if }1/2<x\leq1
\end{array}
\right.
\end{align*}
since the first condition implies $3>1/x\geq2$, i.e., $1>1/x-2\geq0$ and the
second condition implies $2>1/x\geq1$, i.e., $1>1/x-1\geq0$. \ Thus $f$ is
increasing on the left of $x=1/2$ and $f$ is decreasing on the right, which
implies that the median size of $\Lambda_{1}$ is $1/2$. \ 

From
\[%
\begin{array}
[c]{ccc}%
\psi^{\prime}(x)=\dfrac{\varphi(x-1)-\psi(x-1)}{x}, &  & x>2
\end{array}
\]
we have%
\[
\varphi^{\prime}(x)-\psi^{\prime}(x)=-\dfrac{\varphi(x-1)}{x}-\dfrac
{\varphi(x-1)-\psi(x-1)}{x}=\dfrac{-2\varphi(x-1)+\psi(x-1)}{x}%
\]
(a lemmata needed shortly) and%
\[%
\begin{array}
[c]{ccc}%
\psi^{\prime}\left(  \dfrac{1}{x}\right)  =\dfrac{\varphi\left(  \dfrac{1}%
{x}-1\right)  -\psi\left(  \dfrac{1}{x}-1\right)  }{\dfrac{1}{x}}, &  &
0<x<1/2
\end{array}
\]
hence the density $g(x)$ is%
\begin{align*}
\dfrac{d}{dx}\psi\left(  \dfrac{1}{x}\right)   &  =x\left(  \varphi\left(
\dfrac{1}{x}-1\right)  -\psi\left(  \dfrac{1}{x}-1\right)  \right)  \left(
-\dfrac{1}{x^{2}}\right) \\
&  =\dfrac{\psi\left(  \dfrac{1}{x}-1\right)  -\varphi\left(  \dfrac{1}%
{x}-1\right)  }{x}.
\end{align*}
Also, from%
\begin{align*}
\psi^{\prime\prime}(x)  &  =-\dfrac{\varphi(x-1)-\psi(x-1)}{x^{2}}%
+\frac{\varphi^{\prime}(x-1)-\psi^{\prime}(x-1)}{x}\\
&  =\dfrac{-\varphi(x-1)+\psi(x-1)}{x^{2}}+\frac{-2\varphi(x-2)+\psi
(x-2)}{x(x-1)}%
\end{align*}
(by the lemmata) we have
\[
\psi^{\prime\prime}\left(  \dfrac{1}{x}\right)  =\dfrac{-\varphi\left(
\dfrac{1}{x}-1\right)  +\psi\left(  \dfrac{1}{x}-1\right)  }{\dfrac{1}{x^{2}}%
}+\dfrac{-2\varphi\left(  \dfrac{1}{x}-2\right)  +\psi\left(  \dfrac{1}%
{x}-2\right)  }{\dfrac{1}{x}\left(  \dfrac{1}{x}-1\right)  }%
\]
hence (by the chain rule for second derivatives)%
\begin{align*}
\dfrac{d^{2}}{dx^{2}}\psi\left(  \dfrac{1}{x}\right)   &  =\psi^{\prime
}\left(  \dfrac{1}{x}\right)  \dfrac{2}{x^{3}}+\dfrac{1}{x^{4}}\psi
^{\prime\prime}\left(  \dfrac{1}{x}\right) \\
&  =\dfrac{\varphi\left(  \dfrac{1}{x}-1\right)  -\psi\left(  \dfrac{1}%
{x}-1\right)  }{\dfrac{1}{x}}\dfrac{2}{x^{3}}\\
&  +\dfrac{1}{x^{4}}\left[  \dfrac{-\varphi\left(  \dfrac{1}{x}-1\right)
+\psi\left(  \dfrac{1}{x}-1\right)  }{\dfrac{1}{x^{2}}}+\dfrac{-2\varphi
\left(  \dfrac{1}{x}-2\right)  +\psi\left(  \dfrac{1}{x}-2\right)  }{\dfrac
{1}{x}\left(  \dfrac{1}{x}-1\right)  }\right] \\
&  =\dfrac{\varphi\left(  \dfrac{1}{x}-1\right)  -\psi\left(  \dfrac{1}%
{x}-1\right)  }{x^{2}}-\dfrac{2\varphi\left(  \dfrac{1}{x}-2\right)
-\psi\left(  \dfrac{1}{x}-2\right)  }{x^{2}\left(  1-x\right)  }.
\end{align*}
There exists a unique $0<x_{0}<1/2$ for which this expression ($g^{\prime
}(x_{0})$) vanishes. \ Plots of $f(x)$ and $g(x)$ appear in \cite{Fi3-tcs8}
and confirm that $x_{0}$ is the modal size of $\Lambda_{2}$. \ Broadhurst
\cite{Br1-tcs8} obtained an exact equation for $x_{0}$, involving Dickman
dilogarithms and trilogarithms \cite{Br2-tcs8}, then applied numerics. We have
verified his value $x_{0}$ by purely floating point methods.

There is comparatively little to say about medians $\xi_{r}$, defined as
solutions of \cite{KTP-tcs8, DK-tcs8}%
\[
\rho_{r}\left(  \frac{1}{x}\right)  =\frac{1}{2}%
\]
except that $\xi_{1}=1/\sqrt{e}$ is well-known and no closed-form
representations for $\xi_{r}$, $r\geq2$, seem to exist.

\section{Knuth \&\ Trabb Pardo}

An alternative to Heinz's algorithm is one proposed by Knuth \&\ Trabb Pardo
\cite{KTP-tcs8}\ for a restricted case. Define $u_{r}(k,n)$ to be the number
of $n$-permutations whose $r^{\text{th}}$ longest cycle has $\leq k$ nodes
\cite{O0-tcs8}. \ The following recursive formulas apply for $r=1$:
\[
u_{1}(k,n)=\left\{
\begin{array}
[c]{lll}%
{\displaystyle\sum\limits_{m=0}^{k-1}}
\dfrac{(n-1)!}{(n-1-m)!}u_{1}(k,n-1-m) & \bigskip & \text{if }n\geq1\text{ and
}k<n\text{,}\\
n! & \bigskip & \text{if }n\geq1\text{ and }k\geq n\text{,}\\
1 &  & \text{otherwise}%
\end{array}
\right.
\]
and for $r\geq2$:
\[
u_{r}(k,n)=\left\{
\begin{array}
[c]{lll}%
\begin{array}
[c]{l}%
{\displaystyle\sum\limits_{m=0}^{k-1}}
\dfrac{(n-1)!}{(n-1-m)!}u_{r}(k,n-1-m)\;+\\%
{\displaystyle\sum\limits_{m=k}^{n-1}}
\dfrac{(n-1)!}{(n-1-m)!}u_{r-1}(k,n-1-m)
\end{array}
& \bigskip & \text{if }n\geq1\text{ and }k<\left\lfloor n/r\right\rfloor
\text{,}\\
\;n! & \bigskip & \text{if }n\geq1\text{ and }k\geq\left\lfloor
n/r\right\rfloor \text{,}\\
\;1 &  & \text{otherwise.}%
\end{array}
\right.
\]
Clearly $u_{1}(0,n)=\delta_{0,n}$ and $u_{1}(1,n)=1$, hence%
\[
u_{2}(0,4)=u_{1}(0,3)+3u_{1}(0,2)+6u_{1}(0,1)+6u_{1}(0,0)=6.
\]
Also $u_{2}(1,2)=2$ and $u_{2}(1,3)=6$, hence%
\[
u_{2}(1,4)=u_{2}(1,3)+\left[  3u_{1}(1,2)+6u_{1}(1,1)+6u_{1}(1,0)\right]
=6+15=21.
\]
Finally $u_{2}(2,4)=24.$ \ The list%
\[
\left\{  u_{2}(k,4)\right\}  _{k=0}^{2}=\left\{  6,21,24\right\}  =\left\{
6,6+15,21+3\right\}
\]
conveys the same information as the polynomial $p[4,\{0,0\}]$ did in Section
1, although the underlying calculations differed completely.

A\ proof is as follows \cite{KTP-tcs8}. \ We may think of $u_{r}(k,n)$ as
counting permutations on $\{1,\ldots,n\}$ that possess fewer than $r$ cycles
of length exceeding $k$. \ Call such a permutation $(r,n)$\textbf{-good}.
\ Consider now a permutation $P$ on $\{0,1,\ldots,n\}$. \ The node $0$ belongs
to some cycle $C$ within $P$ of length $m+1$. \ Let $P\smallsetminus C$ denote
the permutation which remains upon exclusion of $C$ from $P$. \ Suppose $0\leq
m\leq k-1$; then $P$ is $(r,n+1)$-good iff $P\smallsetminus C$ is
$(r,n-m)$-good. \ Suppose $k\leq m\leq n$; then $P$ is $(r,n+1)$-good iff
$P\smallsetminus C$ is $(r-1,n-m)$-good. \ Thus the formula%
\[
u_{r}(k,n+1)=%
{\displaystyle\sum\limits_{m=0}^{k-1}}
\dfrac{n!}{(n-m)!}u_{r}(k,n-m)\;+%
{\displaystyle\sum\limits_{m=k}^{n}}
\dfrac{n!}{(n-m)!}u_{r-1}(k,n-m)
\]
is true because $n!/(n-m)!$ is the number of possible choices for $C$. \ 

An analog of this recursion for mappings remains open, as far as is known.
\ Finding the number of possible choices for a component $C$ containing the
node $0$ is more complicated than for a cycle containing $0$. \ Each component
consists of a cycle with trees attached; each tree is rooted at a cyclic point
but is otherwise made up of transient points. \ We must account for the
position of $0$ (cyclic or transient?) and the overall configuration
(inventory of tree types and sizes?) \ It would be helpful to learn about
progress in enumerating such $C$ or, if this is impractical, some other
procedure for moving forward.

\section{Une conjecture correspondante}

Short cycles have always presented more analytical difficulties than long
cycles; this paper offers no exception. \ Everything in this section is
conjectural only. \ Define $v_{r}(k,n)$ to be the number of $n$-permutations
whose $r^{\text{th}}$ shortest cycle has $\geq k$ nodes \cite{O0-tcs8}. \ The
following recursive formulas would seem to apply for $r=1$:
\[
v_{1}(k,n)=\left\{
\begin{array}
[c]{lll}%
n! & \bigskip & \text{if }n\geq1\text{ and }k=0\text{,}\\%
{\displaystyle\sum\limits_{m=k-1}^{n-1}}
\dfrac{(n-1)!}{(n-1-m)!}v_{1}(k,n-1-m) & \medskip & \text{if }n\geq1\text{ and
}0<k\leq n\text{,}\\
0 & \medskip & \text{if }n\geq1\text{ and }k>n\text{,}\\
1 &  & \text{otherwise}%
\end{array}
\right.
\]
and for $r\geq2$:
\[
v_{r}(k,n)=\left\{
\begin{array}
[c]{lll}%
\;n! & \bigskip & \text{if }n\geq0\text{ and }k=0\text{,}\\%
\begin{array}
[c]{l}%
\Delta_{r}(k,n)\;+\\%
{\displaystyle\sum\limits_{m=0}^{k-2}}
\dfrac{(n-1)!}{(n-1-m)!}v_{r-1}(k,n-1-m)\;+\\%
{\displaystyle\sum\limits_{m=k-1}^{n-1}}
\dfrac{(n-1)!}{(n-1-m)!}v_{r}(k,n-1-m)
\end{array}
& \medskip & \text{if }n\geq1\text{ and }0<k\leq n-r+1\text{,}\\
\;0 &  & \text{otherwise.}%
\end{array}
\right.
\]
The surprising new term $\Delta_{r}(k,n)$ has a simple formula for $r=2$:
\[%
\begin{array}
[c]{ccccccc}%
\Delta_{2}(k,n)=(n-1)!H_{n-k}, &  & \text{where} &  &
{\displaystyle\sum\limits_{i=1}^{j}}
\dfrac{1}{i}=H_{j}, &  &
{\displaystyle\sum\limits_{i=1}^{j}}
\dfrac{1}{i^{s}}=H_{j,s}%
\end{array}
\]
and unexpected recursions for $r=3$ and $r=4$:
\[
\Delta_{r}(k,n)=\left\{
\begin{array}
[c]{lll}%
\dfrac{1}{2}(n-1)!\left(  H_{n-1}^{2}-H_{n-1,2}\right)  & \bigskip & \text{if
}r=3\text{ and }k=1,\\
\dfrac{1}{6}(n-1)!\left(  H_{n-1}^{3}-3H_{n-1}H_{n-1,2}+2H_{n-1,3}\right)  &
\bigskip & \text{if }r=4\text{ and }k=1,\\
\Delta_{r}(k-1,n)-\dfrac{\Delta_{r-1}(k,n)}{n-k+1} & \medskip & \text{if
}k\geq2\text{ and }n\geq k,\\
0 &  & \text{otherwise.}%
\end{array}
\right.
\]
The values $\Delta_{r}(1,n)$ are unsigned Stirling numbers of the first kind,
i.e., the number of $n$-permutations that have exactly $r$ cycles. \ (Why
should these appear here?)

A\ plausibility argument supporting $v_{r}$ bears resemblance to the proof
underlying $u_{r}$. \ We may think of $v_{r}(k,n)$ as counting permutations on
$\{1,\ldots,n\}$ that possess fewer than $r$ cycles of length surpassed by
$k$. \ Call such a permutation $(r,n)$\textbf{-bad}. \ Let $P$ \& $C$ (of
lengths $n+1$ \& $m+1$) be as before. \ Suppose $0\leq m\leq k-2$; then $P$ is
$(r,n+1)$-bad iff $P\smallsetminus C$ is $(r-1,n-m)$-bad. \ Suppose $k-1\leq
m\leq n$; then $P$ is $(r,n+1)$-bad iff $P\smallsetminus C$ is $(r,n-m)$-bad.
\ This would suggest%
\[
v_{r}(k,n+1)=\Delta_{r}+%
{\displaystyle\sum\limits_{m=0}^{k-2}}
\dfrac{n!}{(n-m)!}v_{r-1}(k,n-m)\;+%
{\displaystyle\sum\limits_{m=k-1}^{n}}
\dfrac{n!}{(n-m)!}v_{r}(k,n-m)
\]
is true with $\Delta_{r}=0$, but experimental data contradict such an
assertion. \ 

Let us illustrate via example, in parallel with Section 3. \ As preliminary
steps, $v_{1}(k,0)=1$ and $v_{1}(n+1,n)=0$, hence%
\[%
\begin{array}
[c]{ccc}%
v_{1}(2,3)=2v_{1}(2,1)+2v_{1}(2,0)=2, &  & v_{1}(3,3)=2v_{1}(3,0)=2.
\end{array}
\]
Clearly $v_{2}(0,4)=24$. \ Also $v_{2}(n,n)=\delta_{0,n}$ and $v_{2}%
(n+1,n)=v_{2}(n+2,n)=0$, hence%
\[
v_{2}(1,2)=\Delta_{2}(1,2)+\left[  v_{2}(1,1)+v_{2}(1,0)\right]  =1+0=1,
\]%
\[
v_{2}(1,3)=\Delta_{2}(1,3)+\left[  v_{2}(1,2)+2v_{2}(1,1)+2v_{2}(1,0)\right]
=3+1=4,
\]%
\[
v_{2}(1,4)=\Delta_{2}(1,4)+\left[  v_{2}(1,3)+3v_{2}(1,2)+6v_{2}%
(1,1)+6v_{2}(1,0)\right]  =11+7=18.
\]
Finally%
\[
v_{2}(2,4)=\Delta_{2}(2,4)+v_{1}(2,3)+\left[  3v_{2}(2,2)+6v_{2}%
(2,1)+6v_{2}(2,0)\right]  =9+2+0=11,
\]%
\[
v_{2}(3,4)=\Delta_{2}(3,4)+\left[  v_{1}(3,3)+3v_{1}(3,2)\right]  +\left[
6v_{2}(3,1)+6v_{2}(3,0)\right]  =6+2+0=8.
\]
Again, the list%
\[
\left\{  v_{2}(k,4)\right\}  _{k=1}^{3}=\left\{  18,11,8\right\}  =\left\{
24-6,18-7,11-3=8\right\}
\]
conveys the same information as the polynomial $q[4,\{\infty,\infty\}]$ did in
Section 1. \ Without the nonzero contribution of $\Delta_{r}(k,n)$, our
modification of Knuth \&\ Trabb Pardo would yield results incompatible with Heinz.

\section{Permutations}

Here \cite{O1-tcs8} are numerical results for $r=2$:

\begin{center}%
\begin{tabular}
[c]{|c|c|c|c|c|c|c|}\hline
$n$ & $_{L}\widetilde{\mu}_{n,2}$ & $_{L}\widetilde{\sigma}_{n,2}^{2}$ &
$_{L}\widetilde{\nu}_{n,2}$ & $_{L}\widetilde{\vartheta}_{n,2}$ &
$_{S}\widetilde{\mu}_{n,2}$ & $_{S}\widetilde{\sigma}_{n,2}^{2}$\\\hline
1000 & 0.209685 & 0.012567 & 0.2110 & 0.2350 & 0.415946 & 1.095918\\\hline
1500 & 0.209650 & 0.012562 & 0.2113 & 0.2353 & 0.408887 & 1.117858\\\hline
2000 & 0.209633 & 0.012560 & 0.2115 & 0.2350 & 0.404309 & 1.131057\\\hline
2500 & 0.209623 & 0.012559 & 0.2112 & 0.2352 & 0.400976 & 1.140134\\\hline
\end{tabular}

Table 5.1:\ Statistics for Permute, rank two ($a=1$)
\end{center}

\noindent as well as $_{S}\nu_{n,2}=2$ for $n>17$ and $_{S}\vartheta_{n,2}=1$
for $n>4$. \ Also%
\[
\lim_{n\rightarrow\infty}\dfrac{_{L}\mu_{n,2}}{n}=\,_{L}G_{1}%
(2,1)=0.20958087428418581398...,
\]%
\[
\lim_{n\rightarrow\infty}\dfrac{_{L}\sigma_{n,2}^{2}}{n^{2}}=\,_{L}%
G_{1}(2,2)-\,_{L}G_{1}(2,1)^{2}=0.01255379063590587814...,
\]%
\[
\lim_{n\rightarrow\infty}\dfrac{_{L}\nu_{n,2}}{n}=\xi_{2}%
=0.21172114641298273896...,
\]%
\[
\lim_{n\rightarrow\infty}\dfrac{_{L}\vartheta_{n,2}}{n}=x_{0}%
=0.23503964593509109370...,
\]%
\[
\lim_{n\rightarrow\infty}\dfrac{_{S}\mu_{n,2}}{\ln(n)^{2}}=\frac{e^{-\gamma}%
}{2}=0.28072974178344258491...,
\]%
\[
\lim_{n\rightarrow\infty}\dfrac{_{S}\sigma_{n,2}^{2}}{n\ln(n)}=\,_{S}%
G_{P}(2,2)=1.30720779891056809974....
\]

\noindent The final $n\ln(n)$ asymptotic is based on \cite{SL-tcs8, Shi-tcs8},
not (inaccurate) Theorem 5 in \cite{PR2-tcs8}.

Here \cite{O2-tcs8} are numerical results for $r=3$:

\begin{center}%
\begin{tabular}
[c]{|c|c|c|c|c|c|c|}\hline
$n$ & $_{L}\widetilde{\mu}_{n,3}$ & $_{L}\widetilde{\sigma}_{n,3}^{2}$ &
$_{L}\widetilde{\nu}_{n,3}$ & $_{L}\widetilde{\vartheta}_{n,3}$ &
$_{S}\widetilde{\mu}_{n,3}$ & $_{S}\widetilde{\sigma}_{n,3}^{2}$\\\hline
1000 & 0.088357 & 0.004499 & 0.0750 & 0.0010 & 0.155997 & 0.450101\\\hline
1500 & 0.088344 & 0.004497 & 0.0753 & 0.0007 & 0.153079 & 0.468681\\\hline
2000 & 0.088337 & 0.004496 & 0.0755 & 0.0005 & 0.151161 & 0.480325\\\hline
2500 & 0.088333 & 0.004496 & 0.0756 & 0.0004 & 0.149752 & 0.488548\\\hline
\end{tabular}

Table 5.2:\ Statistics for Permute, rank three ($a=1$)
\end{center}

\noindent as well as $_{S}\nu_{n,2}=7$ for $n>370$ and $_{S}\vartheta_{n,2}=2$
for $n>49$. \ Also%
\[
\lim_{n\rightarrow\infty}\dfrac{_{L}\mu_{n,3}}{n}=\,_{L}G_{1}%
(3,1)=0.08831609888315363101...,
\]%
\[
\lim_{n\rightarrow\infty}\dfrac{_{L}\sigma_{n,3}^{2}}{n^{2}}=\,_{L}%
G_{1}(3,2)-\,_{L}G_{1}(3,1)^{2}=0.00449392318179080474...,
\]%
\[
\lim_{n\rightarrow\infty}\dfrac{_{L}\nu_{n,3}}{n}=\xi_{3}%
=0.07584372316630152789...,
\]%
\[
\lim_{n\rightarrow\infty}\dfrac{_{L}\vartheta_{n,3}}{n}=0,
\]%
\[
\lim_{n\rightarrow\infty}\dfrac{_{S}\mu_{n,3}}{\ln(n)^{3}}=\frac{e^{-\gamma}%
}{6}=0.09357658059448086163...,
\]%
\[
\lim_{n\rightarrow\infty}\dfrac{_{S}\sigma_{n,3}^{2}}{n\ln(n)^{2}}=\,_{S}%
G_{P}(3,2)=0.65360389945528404987....
\]

\noindent The final $n\ln(n)^{2}$ asymptotic is based on \cite{SL-tcs8,
Shi-tcs8}.

Here \cite{O3-tcs8} are numerical results for $r=4$:

\begin{center}%
\begin{tabular}
[c]{|c|c|c|c|c|c|c|}\hline
$n$ & $_{L}\widetilde{\mu}_{n,4}$ & $_{L}\widetilde{\sigma}_{n,4}^{2}$ &
$_{L}\widetilde{\nu}_{n,4}$ & $_{L}\widetilde{\vartheta}_{n,4}$ &
$_{S}\widetilde{\mu}_{n,4}$ & $_{S}\widetilde{\sigma}_{n,4}^{2}$\\\hline
1000 & 0.040353 & 0.001586 & 0.0260 & 0.0010 & 0.042215 & 0.118491\\\hline
1500 & 0.040351 & 0.001585 & 0.0267 & 0.0007 & 0.041482 & 0.126180\\\hline
2000 & 0.040349 & 0.001585 & 0.0265 & 0.0005 & 0.040987 & 0.131244\\\hline
2500 & 0.040348 & 0.001585 & 0.0268 & 0.0004 & 0.040618 & 0.134938\\\hline
\end{tabular}

Table 5.3:\ Statistics for Permute, rank four ($a=1$)
\end{center}

\noindent as well as $_{S}\nu_{n,4}=19$ for $n>1482$ and $_{S}\vartheta
_{n,4}=3$ for $n>666$. \ Also%
\[
\lim_{n\rightarrow\infty}\dfrac{_{L}\mu_{n,4}}{n}=\,_{L}G_{1}%
(4,1)=0.04034198873687046287...,
\]%
\[
\lim_{n\rightarrow\infty}\dfrac{_{L}\sigma_{n,4}^{2}}{n^{2}}=\,_{L}%
G_{1}(4,2)-\,_{L}G_{1}(4,1)^{2}=0.00158383677354017280...,
\]%
\[
\lim_{n\rightarrow\infty}\dfrac{_{L}\nu_{n,4}}{n}=\xi_{4}%
=0.02713839684981404992...,
\]%
\[
\lim_{n\rightarrow\infty}\dfrac{_{L}\vartheta_{n,4}}{n}=0,
\]%
\[
\lim_{n\rightarrow\infty}\dfrac{_{S}\mu_{n,4}}{\ln(n)^{4}}=\frac{e^{-\gamma}%
}{24}=0.02339414514862021540...,
\]%
\[
\lim_{n\rightarrow\infty}\dfrac{_{S}\sigma_{n,4}^{2}}{n\ln(n)^{3}}=\,_{S}%
G_{P}(4,2)=0.21786796648509468329....
\]
The final $n\ln(n)^{3}$ asymptotic is based on \cite{SL-tcs8, Shi-tcs8}.

\section{Mappings}

Our modified Knuth \&\ Trabb Pardo algorithm is unavailable in this setting,
thus we turn to Heinz's program. \ A general observation for $2\leq r\leq4$ is
$_{L}\vartheta_{n,r}=0$ always. \ Here \cite{O4-tcs8} are numerical results
for $r=2$:

\begin{center}%
\begin{tabular}
[c]{|c|c|c|c|c|c|c|}\hline
$n$ & $_{L}\widetilde{\mu}_{n,2}$ & $_{L}\widetilde{\sigma}_{n,2}^{2}$ &
$_{L}\widetilde{\nu}_{n,2}$ & $_{S}\widetilde{\mu}_{n,2}$ & $_{S}%
\widetilde{\sigma}_{n,2}^{2}$ & $_{S}\widetilde{\nu}_{n,2}$\\\hline
100 & 0.166817 & 0.019535 & 0.1300 & 0.680589 & 0.279032 & 0.1200\\\hline
200 & 0.168100 & 0.019243 & 0.1400 & 0.718071 & 0.323910 & 0.0750\\\hline
300 & 0.168642 & 0.019121 & 0.1433 & 0.737331 & 0.350358 & 0.0567\\\hline
400 & 0.168959 & 0.019050 & 0.1450 & 0.749928 & 0.368810 & 0.0450\\\hline
\end{tabular}

Table 6.1:\ Statistics for Map, rank two ($a=1/2$)
\end{center}

\noindent as well as $_{S}\nu_{n,2}=19$ for $n>443$ and $_{S}\vartheta
_{n,2}=2$ for $n>433$. \ Let us elaborate on the latter statistic (because it
seems surprising at first glance:\ an extended string of $0$s abruptly
switches to $2$s). \ If $\pi_{r}(k,n)$ denotes the probability that the
$r^{\text{th}}$ smallest component of a random $n$-mapping has exactly $k$
nodes, then%
\[%
\begin{array}
[c]{l}%
\left\{  \pi_{2}(k,432\right\}  _{k=0}^{4}=\left\{
0.0595400,0.0532617,0.0594378,0.0477544,0.0387585\right\}  ,\\
\left\{  \pi_{2}(k,433\right\}  _{k=0}^{4}=\left\{
0.0594720,0.0532614,0.0594373,0.0477539,0.0387581\right\}  ,\\
\left\{  \pi_{2}(k,434\right\}  _{k=0}^{4}=\left\{
0.0594044,0.0532612,0.0594369,0.0477535,0.0387576\right\}  ,\\
\left\{  \pi_{2}(k,435\right\}  _{k=0}^{4}=\left\{
0.0593369,0.0532609,0.0594365,0.0477530,0.0387571\right\}  .
\end{array}
\]
The maximum probability clearly is at $k=0$ for $n\leq433$ and then shifts to
$k=2$ for $n\geq434$. \ Also
\[
\lim_{n\rightarrow\infty}\dfrac{_{L}\mu_{n,2}}{n}=\,_{L}G_{1/2}%
(2,1)=0.17090961985966239214...,
\]%
\[
\lim_{n\rightarrow\infty}\dfrac{_{L}\sigma_{n,2}^{2}}{n^{2}}=\,_{L}%
G_{1/2}(2,2)-\,_{L}G_{1/2}(2,1)^{2}=0.01862022330678138872...,
\]%
\[
\lim_{n\rightarrow\infty}\dfrac{_{L}\nu_{n,2}}{n}=0.148...,
\]%
\[
\lim_{n\rightarrow\infty}\dfrac{_{S}\mu_{n,2}}{n^{1/2}\ln(n)}=\sqrt{2}%
\,_{S}G_{1/2}(2,1)=2.06089224152016653900...,
\]%
\[
\lim_{n\rightarrow\infty}\dfrac{_{S}\sigma_{n,2}^{2}}{n^{3/2}\ln(n)}=\sqrt
{2}\,_{S}G_{1/2}(2,2)=1.40007638550124502818....
\]
No exact equation (akin to one involving $\rho_{r}$ in Section 2) is known for
the median of $L$. \ An $r^{\text{th}}$ order Dickman function $\rho_{r,1/2}$
of type $a=1/2$ might be needed. \ What is responsible for mismatches between
data and theory for $S$? \ This may be due to uncertainty about how the
correction factor $\sqrt{2}$ should be generalized from $r=1$ to all $r\geq1$.
\ We believe that the sequence $_{S}\nu_{n,2}$ is bounded; a proof is not known.

Here \cite{O5-tcs8} are numerical results for $r=3$: \ 

\begin{center}%
\begin{tabular}
[c]{|c|c|c|c|c|c|c|}\hline
$n$ & $_{L}\widetilde{\mu}_{n,3}$ & $_{L}\widetilde{\sigma}_{n,3}^{2}$ &
$_{L}\widetilde{\nu}_{n,3}$ & $_{S}\widetilde{\mu}_{n,3}$ & $_{S}%
\widetilde{\sigma}_{n,3}^{2}$ & $_{S}\widetilde{\nu}_{n,3}$\\\hline
100 & 0.044147 & 0.003902 & 0 & 0.126620 & 0.052261 & 0.0700\\\hline
150 & 0.045094 & 0.003902 & 0.0067 & 0.133605 & 0.055079 & 0.0867\\\hline
200 & 0.045642 & 0.003903 & 0.0100 & 0.138200 & 0.057284 & 0.0850\\\hline
250 & 0.046008 & 0.003904 & 0.0120 & 0.141572 & 0.059120 & 0.0880\\\hline
\end{tabular}

Table 6.2:\ Statistics for Map, rank three ($a=1/2$)
\end{center}

\noindent as well as $_{S}\nu_{n,3}=24$ for $n>275$ and $_{S}\vartheta
_{n,3}=0$ for $n\leq278$ at least. \ Also%
\[
\lim_{n\rightarrow\infty}\dfrac{_{L}\mu_{n,3}}{n}=\,_{L}G_{1/2}%
(3,1)=0.04889742536845958914...,
\]%
\[
\lim_{n\rightarrow\infty}\dfrac{_{L}\sigma_{n,3}^{2}}{n^{2}}=\,_{L}%
G_{1/2}(3,2)-\,_{L}G_{1/2}(3,1)^{2}=0.00392148747204257695...,
\]%
\[
\lim_{n\rightarrow\infty}\dfrac{_{S}\mu_{n,3}}{n^{1/2}\ln(n)^{2}}=\sqrt
{2}\,_{S}G_{1/2}(3,1)=1.03044612076008326950...,
\]%
\[
\lim_{n\rightarrow\infty}\dfrac{_{S}\sigma_{n,3}^{2}}{n^{3/2}\ln(n)^{2}}%
=\sqrt{2}\,_{S}G_{1/2}(3,2)=0.70003819275062251409....
\]
The median of $L$ is unknown and mismatches worsen. It is certainly possible
that the sequence $_{L}\nu_{n,3}$ might be bounded; the trend of $_{S}%
\nu_{n,3}$ is ambiguous. \ There are presently insufficient data to render judgement.

Here \cite{O6-tcs8} are numerical results for $r=4$:

\begin{center}%
\begin{tabular}
[c]{|c|c|c|c|c|}\hline
$n$ & $_{L}\widetilde{\mu}_{n,4}$ & $_{L}\widetilde{\sigma}_{n,4}^{2}$ &
$_{S}\widetilde{\mu}_{n,4}$ & $_{S}\widetilde{\sigma}_{n,4}^{2}$\\\hline
100 & 0.011968 & 0.000710 & 0.015300 & 0.007424\\\hline
125 & 0.012324 & 0.000717 & 0.016032 & 0.007682\\\hline
150 & 0.012585 & 0.000722 & 0.016606 & 0.007877\\\hline
175 & 0.012787 & 0.000726 & 0.017077 & 0.008034\\\hline
\end{tabular}

Table 6.3:\ Statistics for Map, rank four ($a=1/2$)
\end{center}

\noindent as well as $_{L}\nu_{n,4}=0$, $_{S}\nu_{n,4}=0$, $_{S}%
\vartheta_{n,4}=0$ for $n\leq183$ at least.. \ Also%
\[
\lim_{n\rightarrow\infty}\dfrac{_{L}\mu_{n,4}}{n}=\,_{L}G_{1/2}%
(4,1)=0.01514572139988693564...,
\]%
\[
\lim_{n\rightarrow\infty}\dfrac{_{L}\sigma_{n,4}^{2}}{n^{2}}=\,_{L}%
G_{1/2}(4,2)-\,_{L}G_{1/2}(4,1)^{2}=0.00077636923173854484...,
\]%
\[
\lim_{n\rightarrow\infty}\dfrac{_{S}\mu_{n,4}}{n^{1/2}\ln(n)^{3}}=\sqrt
{2}\,_{S}G_{1/2}(4,1)=0.34348204025336108983...,
\]%
\[
\lim_{n\rightarrow\infty}\dfrac{_{S}\sigma_{n,4}^{2}}{n^{3/2}\ln(n)^{3}}%
=\sqrt{2}\,_{S}G_{1/2}(4,2)=0.23334606425020750469....
\]
Again, the median of $L$ is unknown and mismatches worsen. \ Although both
sequences $_{L}\nu_{n,4}$ and $_{S}\nu_{n,4}$ seem to be bounded (only $0$s
observed), we sense that they are still in transience and substantially more
data will be required to reach steady state.

\section{Acknowledgements}

I am indebted to Alois Heinz for providing the algorithms underlying
$p[n,\ell]$ and $q[n,\ell]$, and to David Broadhurst for calculating
$_{L}\widetilde{\nu}_{n,2}$, $_{L}\widetilde{\nu}_{n,3}$, $_{L}\widetilde{\nu
}_{n,4}$, $_{L}\widetilde{\vartheta}_{n,2}$ to high precision as
$n\rightarrow\infty$ (permutations only). \ Many thanks are owed to
Jean-Francois Alcover for translating Heinz's concise Maple code to a form
I\ could understand. The volunteers who edit and maintain OEIS, the creators
of Mathematica, as well as administrators of the MIT Engaging Cluster, earn my
gratitude every day. \ A sequel to this paper will be released soon
\cite{Fi4-tcs8}.

\end{document}